\documentclass{article}

\begin{document}

\newtheorem{thm}{Theorem}
\newtheorem{prop}{Proposition}
\newtheorem{defn}{Definition}
\newtheorem{condn}{Condition}
\newtheorem{cor}{Corollary}
\newtheorem{obs}{Observation}
\newcommand{\tu}{\tilde{u}}
\newcommand{\ux}{{\underline{x}}}
\newcommand{\utu}{\underline{\tilde{u}}}
\newcommand{\uh}{\underline{h}}
\newcommand{\oh}{\overline{h}}
\newcommand{\oP}{\overline{P}}
\newcommand{\tP}{\tilde{P}}
\newcommand{\tA}{\tilde{A}}
\newcommand{\ut}{\underline{t}}
\newcommand{\ot}{\overline{t}}
\newcommand{\tsigma}{\tilde{\sigma}}
\newcommand{\tU}{\tilde{U}}
\newcommand{\ufl}{u_{\flat}}
\newcommand{\xifl}{\xi^{\flat}}
\newcommand{\ua}{\underline{a}}
\newcommand{\Pfl}{P_{\flat}}
\newcommand{\Afl}{A_{\flat}}
\newcommand{\Vfl}{V_{\flat}}
\newcommand{\Area}{{\rm Area}\ }
\newcommand{\oM}{\overline{M}}
\newcommand{\bR}{{\bf R}}
\newcommand{\bC}{{\bf C}}
\newcommand{\uu}{\underline{u}}
\newcommand{\Vol}{{\rm Vol}}
\newcommand{\Mmax}{M_{{\rm max}}}
\newcommand{\Fmax}{F_{{\rm max}}}
\newcommand{\udelta}{\underline{\delta}}
\newcommand{\xist}{\xi^{*}}
\newcommand{\distu}{{\rm dist}_{u}}
\newcommand{\distual}{{\rm dist}_{u^{(\alpha)}}}
\newcommand{\uZ}{\underline{Z}}
\newcommand{\uF}{\underline{F}}

\title{A generalised Joyce construction for a family of nonlinear partial
differential equations}
\author{S. K. Donaldson}
\date{\today}
\maketitle

The construction of Dominic Joyce referred to in the title is for
 self-dual Riemannian metrics with two Killing fields \cite{kn:J}.
In this note we show that this construction (in the form described by Calderbank
and Pedersen in \cite{kn:CP}) can be extended to a family of
nonlinear fourth order PDE in two dimensions, essentially reducing them to
linear equations. One equation in this family is the {\it affine maximal equation},
and we will see that the construction in this case is amost the same as one due to Chern and Terng.
\section{The main result}

We consider convex functions $u(x_{i})$ defined on a domain in $\bR^{n}$
and write 
$J=\det(u_{ij})$, where $(u_{ij})$ is the Hessian $\left( \frac{\partial^{2}u}{\partial
x_{i}\partial x_{j}}\right)$. Let $\psi$ be any smooth, strictly convex function on the half-line
$(0,\infty)$ and consider the functional
\begin{equation} {\cal F}= \int \psi(J)\ dx_{1}\dots dx_{n}. \end{equation}
The corresponding Euler-Lagrange equations $\delta {\cal F}=0$ are
\begin{equation} \sum_{ij} \frac{\partial^{2}}{\partial x_{i}\partial x_{j}}\left( J \psi'(J) u^{ij}\right) =0, \end{equation}
where $(u^{ij})$ is the matrix inverse of $(u_{ij})$. This a  nonlinear
fourth order PDE for the function $u$. In the case when $\phi(J)= -J^{\alpha}$
for some $\alpha\in (0,1)$ these equations have been studied by Trudinger
and Wang \cite{kn:TW}. The case covered by Joyce's original construction is when $n=2$ and $\psi(J)=-\log J$, as we will discuss further in Section
3. (Of course we only
consider the functional ${\cal F}$ as a motivation for writing down the partial
differential equations (2), and the actual convergence of the integral (1)
is irrelevant.) Let us say that a point $(x_{1}, x_{2})$ is an \lq\lq ordinary
point''  if the derivative $\nabla J$ does not vanish there.

  To state our result, let $f$ be a solution of the equation
  $$   f'(t)= t^{1/2} \psi''(t). $$
  Thus $f'>0$ and $f$ is a $1-1$ map from $(0,\infty)$ to some finite or
  infinite interval $I$. Let $p$ be   the inverse map raised to the power
  $-1/2$, so
  $     J^{-1/2} = p(r) $ when $r=f(J)$. Now consider the linear, second order
  PDE for a function $\xi(H,r)$ defined on some domain in $\bR\times I$;
  \begin{equation}\
  \frac{\partial^{2} \xi}{\partial H^{2}} + \frac{1}{p(r)} \frac{\partial
  }{\partial r} \left( p(r) \frac{\partial \xi}{\partial r}\right)= 0.\end{equation}
  
  \begin{thm}
  Suppose $\xi_{1}, \xi_{2}$ are two solutions of equation (3) and $\det\left(\frac{\partial(\xi_{1},\partial
  \xi_{2})}{\partial(H,r)}\right)$ is positive at some point $( H_{0}, r_{0})$. Define
  three $1$-forms $ \epsilon_{1}, \epsilon_{2}, \epsilon$ by
  $$  \epsilon_{1} = p(r)( \frac{\partial \xi_{2}}{\partial r} dH- \frac{\partial
  \xi_{2}}{\partial H} dr) \ \ \ \epsilon_{2}= p(r)(- \frac{\partial \xi_{1}}{\partial
  r} dH + \frac{\partial \xi_{1}}{\partial H} dr)$$
  $$  \epsilon= \xi_{1} \epsilon_{1} + \xi_{2} \epsilon_{2}. $$
  Then 
   $\epsilon, \epsilon_{1}, \epsilon_{2}$ are closed $1$-forms and $\epsilon_{1}\wedge\epsilon_{2}$
  is non-vanishing near $(H_{0}, r_{0})$. Thus  we can find functions $u, x_{1},
  x_{2}$ with
  $du=\epsilon, dx_{1}=\epsilon_{1}, dx_{2}=\epsilon_{2}$ and
  $x_{1}, x_{2}$ give local co-ordinates around $(H_{0}, r_{0})$. If
  we regard $u$ as a function of $(x_{1}, x_{2})$ it is a solution of the
  equation (2), on a suitable domain. All points in this domain are ordinary
  points. Conversely,   any solution of (2) with $n=2$,   in the neighbourhood of an ordinary
  point, is obtained in this way, with solutions $\xi_{1}, \xi_{2}$ of (3) which
  are unique up to translation in the $H$-variable and the addition of constants
  \end{thm}
  
  In sum, the local study of the nonlinear equation (2), in  two dimensions, is essentially equivalent
  to that of the linear equation (3).
  
  \section{The proof}
  
  This is entirely elementary. From now on we always suppose the dimension
  $n$ is $2$, and we work locally so we will not specify the precise domains
  of definition of the various functions. Given a convex function $u$ we consider the Riemannian metric
  $g=\sum u_{ij} dx_{i} dx_{j}$, and in particular the {\it conformal structure}
  which this defines. Write $\xi_{i}$ for the Legendre transform
  co-ordinates $\xi_{i}= \frac{\partial u}{\partial x_{i}}$. Suppose now
  that
  we have some other local co-ordinates $\lambda_{1}, \lambda_{2}$, inducing
  the same orientation as $x_{1}, x_{2}$ (i.e. $\det\left(\frac{\partial x_{i}}{\partial
  \lambda_{a}}\right)>0$). We
  can write the metric $g=\sum g_{ab} d\lambda_{a}d\lambda_{b}$ in these
  co-ordinates. Recall that the co-ordinates are called {\it isothermal}
  if the matrix $(g_{ab})$ is a multiple of the identity matrix at each point,
  or in other words
  $$   g= V (d\lambda_{1}^{2} +d\lambda_{2}^{2}), $$
  for a positive function $V(\lambda_{1}, \lambda_{2})$. Let $\epsilon_{ij}$
  denote the alternating tensor  $\epsilon_{12}= -\epsilon_{21}=1,
\epsilon_{11}=\epsilon_{22}=0$.
 
  \begin{obs}
  The co-ordinates $\lambda_{1}, \lambda_{2}$ are isothermal if and only
  if the partial derivatives $\frac{\partial x_{i}}{\partial \lambda_{a}},
  \frac{\partial \xi_{j}}{\partial \lambda_{b}}$ are related by the four equations
  \begin{equation}
  \frac{\partial \xi_{i}}{\partial \lambda_{a}}=  \epsilon_{a b}\epsilon_{i j} \sqrt{J}\frac{\partial
x_{j}}{\partial \lambda_{b}}, \end{equation}
where $J=J(\ux(\lambda_{1}, \lambda_{2}))$.

\end{obs}

To see this we write
$$  g=\sum_{ij} u_{ij} dx_{i} dx_{j} = \sum_{j} d\xi_{j} dx_{j} =\sum_{j,a,b} \frac{\partial \xi_{j}}{\partial \lambda_{a}} \frac{\partial x_{j}}{\partial
\lambda_{b}} d\lambda_{a} d\lambda_{b}. $$
So the isothermal condition is
$$   \sum_{j} \frac{\partial \xi_{j}}{\partial \lambda_{a}} \frac{\partial
x_{j}}{\partial \lambda_{b}}= V \delta_{ab}. $$
In matrix notation, if $A=\left(\frac{\partial x_{i}}{\partial \lambda_{a}}\right),
B=\left(\frac{\partial \xi_{i}}{\partial \lambda_{a}}\right)$, this is
$$    A^{T} B = V 1, $$
or in other words
\begin{equation}  B^{T}  = V  A^{-1}. \end{equation}
Taking determinants we have $\det B \det A=V^{2}$. On the other hand the
matrix $(u_{ij})$ is $BA^{-1}$, so $J=\det B \det A^{-1}$. Thus $V=\sqrt{J}
\det A$ and (5) is
$$  B^{T}= \sqrt{J}  \det A\  A^{-1}, $$
which is the same as (4), by the formula for the inverse of a $2\times 2$ matrix.

\

\

Now consider any positive function $P(\lambda_{1}, \lambda_{2})$ and the
pair of second order linear PDE
\begin{equation}  \frac{\partial}{\partial \lambda_{1}}\left( P \frac{\partial\xi}
 {\partial \lambda_{1}}\right) + \frac{\partial}{\partial \lambda_{2}}\left(
 P
 \frac{\partial \xi}{\partial \lambda_{2}}\right) = 0, \end{equation}
 
\begin{equation}  \frac{\partial}{\partial \lambda_{1}}\left( P^{-1} \frac{\partial
x}
 {\partial \lambda_{1}}\right) + \frac{\partial}{\partial \lambda_{2}}\left(
 P^{-1}
 \frac{\partial x }{\partial \lambda_{2}}\right) = 0. \end{equation}

 \begin{obs}
 If $\xi(\lambda_{1},\lambda_{2})$ is a solution of (6) then there is a solution
 $x(\lambda_{1}, \lambda_{2})$ of the first order system
 \begin{equation}  \frac{\partial x}{\partial \lambda_{a}} = P \epsilon_{ab} \frac{\partial \xi}{\partial \lambda_{b}} , \end{equation}
unique up to the addition of a constant, and $x(\lambda_{1},\lambda_{2})$
satisfies (7).
\end{obs}
This is because the consistency condition for the first order system is
$$  \frac{\partial }{\partial \lambda_{2}} \left( P\frac{\partial \xi}{\partial
\lambda_{2}}\right)= \frac{\partial}{\partial \lambda_{1}}\left(-P \frac{\partial \xi}{\partial \lambda_{1}}\right), $$
which is the equation (6).  In a different language, we  are saying that the $1$-form
$\epsilon= P^{-1} \frac{\partial \xi}{\partial \lambda}_{2} d\lambda_{1}-
P \frac{\partial \xi}{\partial \lambda_{1}}d\lambda_{2}$ is closed, so can be written
as $dx$, for a function $x$. Note that, changing $P$ to $P^{-1}$, there is a complete
symmetry between $x$ and $\xi$ so we can also start with a solution of (7)
and construct a solution of (6).

\

\begin{obs}
Suppose $\xi_{1}, \xi_{2}$ are two solutions of (6), with $\det\left(\frac{\partial \xi_{i}}{\partial
\lambda_{a}}\right) >0$. Let $x_{2}$ be the solution of (7) corresponding to $\xi_{1}$
by (8), and $x_{1}$ be the solution corresponding to $-\xi_{2}$. Then
$\det\left( \frac{\partial x_{i}}{\partial \lambda_{a}}\right)>0$, so $x_{1},
x_{2}$ give local co-ordinates. Write
$$   \epsilon= \xi_{1} dx_{1} + \xi_{2} dx_{2}.$$
Then $\epsilon$ is a closed $1$-form and so $\epsilon =du$ for a function
$u$. If we express $u$ as a function of $x_{1}, x_{2}$ then
$\frac{\partial u}{\partial x_{i}}= \xi_{i}$. 
\end{obs}

This is all straightforward. The conditions (8) imply that $$\det\left(\frac{\partial
x_{i}}{\partial \lambda_{a}}\right)=P^{2}\det\left( \frac{\partial \xi_{i}}{\partial
\lambda_{a}}\right)>0.$$
We have
$  d\epsilon= d\xi_{1}dx_{1}+d\xi_{2} dx_{2}  $ and, 
writing $\xi_{i,a}=\frac{\partial \xi_{i}}{\partial \lambda_{a}}$, 
$$  d\epsilon= (\xi_{1,1} d\lambda_{1}
+ \xi_{1,2} d\lambda_{2})(-\xi_{2,2} d\lambda_{1}+\xi_{2,1} d\lambda_{2})+
(\xi_{2,1} d\lambda_{1} +\xi_{2,2} d\lambda_{2})(\xi_{1,2} d\lambda_{1}-
\xi_{1,1} d\lambda_{2}) = 0. $$

\

\

Now return to our functions $\psi(J), f(J)$ and the Euler-Lagrange equation
(2).
\begin{obs}
A convex function $u$ satisfies equation (2) if and only if $f(J)$ is harmonic
with respect to the metric $g=\sum u_{ij} dx_{i} dx_{j}$.
\end{obs}
This is true in any dimension.
The formula for the derivative of an inverse matrix is
$$  \frac{\partial}{\partial x_{k}} u^{ij}= - \sum_{ p q} u^{ip} u_{pq k}
u_{qj}, $$ whereas the formula for the derivative of the determinant is
$$   \frac{\partial J}{\partial x_{i}}= J \sum_{pq} u^{pq} u_{pq i}. $$
These yield the identity 
\begin{equation} \sum_{j} \frac{\partial}{\partial  x_{j}} u^{ij}= - \sum_{p q j} u^{ip}
u_{j pq} u^{ qj}= -\sum_{p} u^{ip} J^{-1} \frac{\partial J}{\partial x_{p}}.
\end{equation}
Our Euler-Lagrange equation (2) is
\begin{equation}   \sum_{i} \frac{\partial v_{i}}{\partial x_{i}} = 0, \end{equation}
where $$ v_{i}= \sum_{j} \frac{\partial}{\partial x_{j} } ( J \psi'(J)
u^{ij}). $$
So $$v_{i}= \sum_{j} \frac{\partial (J \psi'(J))}{\partial x_{j}} u^{ij} + J
\psi'(J) \frac{\partial u^{ij}}{\partial x_{j}}. $$
 By the definition of the function $f$ we have
 $$  \frac{\partial}{\partial x_{j}}( J \psi'(J)) = \sqrt{J} \frac{\partial
 f(J)}{\partial x_{j}} + \psi'(J) \frac{\partial J}{\partial x_{j}}.
 $$
 Using (9) we obtain
 $$  v_{i}= \sqrt{J} \frac{\partial f(J)}{\partial x_{j}}u^{ij}. $$
 Thus the equation (10) is the Laplace equation in the metric $g$:
 $$  \sum_{i} \frac{\partial}{\partial x_{i}}\left( \sqrt{J} u^{ij} \frac{\partial
 f(J)}{\partial x_{j}} \right) = 0. $$

 \

 With these four observations the main result, Theorem 1, is almost obvious. Suppose we start with an ordinary point of a  solution $u$ to (2). Then $r=f(J)$ is harmonic, by Observation 4, 
and we can suppose that its derivative does not vanish in the region considered. There is then a conjugate
harmonic function $H$, which by definition is one such that the local co-ordinates
$(H,r)$ are isothermal. 
By Observations 1 and 2 the functions $x_{i}, \xi_{j}$ satisfy the equations (6),(7)
respectively with $\lambda_{1}=H, \lambda_{2}=r$ and $P=p(r)$. Since $P$
does not depend on $H$ the equation (6) can be written in the form (3). By Observations
2 and 3 we can recover the original function $u$ from the two solutions $\xi_{1},\xi_{2}$
of the linear PDE (or, equally well, the two solutions $x_{1}, x_{2}$). It
is also clear that, conversely, starting with any two solutions $\xi_{1},\xi_{2}$ to the
linear equation we construct a solution to (2) by this method.

\

Note that the conjugate function $H$, in this situation, can be  defined
simply
as the solution of the system
\begin{equation}   \frac{\partial H}{\partial x_{i}} = \epsilon_{ij} v_{j},
\end{equation} so we can also thing of $H$ as the Hamiltonian generating
the area-preserving vector field $v$ in the $(x_{1}, x_{2})$ plane. We leave
this as an exercise for the reader.
\section{Examples and discussion}

\begin{enumerate}
\item  Let $u^{*}(\xi_{1},\xi_{2})$ be the Legendre transform of $u(x_{1},
x_{2})$. The Hessian $(\frac{\partial^{2}u^{*}}{\partial \xi_{i}\partial
\xi_{j}})$ at the point $\xi_{i}=\frac{\partial u}{\partial x_{i}}$ is the
inverse of the Hessian of $u$ at $x_{i}$. Thus we can write
$${\cal F}= \int \psi(J_{*}^{-1}) J_{*}\ d\xi_{1} d\xi_{2}, $$
where $J_{*}=\det(\frac{\partial^{2} u^{*}}{\partial \xi_{i}\partial \xi_{j}})$.
This implies that the Legendre transform takes a solution $u$ of the equation (2) associated
with $\psi$ to a solution $u^{*}$ of the equation associated to the function
$$\psi^{*}(t)= t \psi(t^{-1}). $$
In our construction this just means replacing the function $p(r)$ by $p(-r)^{-1}$
and interchanging the roles of the co-ordinates $x_{i}, \xi_{i}$.

\item If $\psi(t)=-\log t$ the equation (2) is $\sum \frac{\partial^{2}u^{ij}}{\partial
x_{i}\partial x_{j}}=0$. This is the equation defining a zero scalar curvature
Kahler metric in \lq\lq symplectic'' co-ordinates, see \cite{kn:A}, \cite{kn:D}.
Explicitly we introduce two further co-ordinates $\theta_{1}, \theta_{2}$
and consider the Riemannian metric, in four dimensions, 
$$  \sum u_{ij} dx_{i} dx_{j}+ \sum u^{ij} d\theta_{i}d\theta_{j}. $$
 The well-known
fact that, in two complex dimensions,
such metrics are self-dual gives the link with Joyce's original formulation
of his construction. Under the Legendre transform we get another description
corresponding to the function $\psi^{*}(t)=t\log t$. This leads to the equations
describing zero scalar curvature metrics in \lq\lq complex'' co-ordinates.
When $\psi(t)=-\log t$ we get $p(r)=r$ and equation (3) is the familiar equation defining axi-symmetric harmonic functions on $\bR^{3}$.

\item For any function $p(r)$ there is an obvious  solution $\xi=H$ to (3).
Thus we get a special family of solutions to (2) with $\xi_{1}=H$ and $\xi_{2}$
some other solution of (3). There is another special family, which corresponds
to this under the Legendre transform, when $\xi_{1}$ is a function of $r$
only, so $x_{2}=H$. These correspond to second order equations of Monge-Amp\`ere
type
which give special solutions of (2). For example, in the zero scalar curvature
case above we have the special solutions where $\xi_{1}=\log r$. The function
$u$ satisfies the equation $\det \left( \frac{\partial^{2}u}{\partial x_{i}\partial x_{j}}\right)=e^{-\xi_{1}/2}$
and the corresponding four dimensional Riemannian metric is Ricci-flat.

\item For the Trudinger-Wang equations, when $\psi(t)=-t^{\alpha}$ with $0<\alpha<1/2$, we get $p(r)= r^{1/(1-2\alpha)}$ (up to factor which is irrelevant because it drops out
of the equation (3)). If $\alpha=1/2$ we get $p(r)=e^{r/2}$.

\item Consider the graph of a function $u(x_{1}, x_{2})$ as a surface in
$\bR^{3}$. The Gauss curvature is
$$  K= \frac{\det(u_{ij})}{(1+\vert \nabla u \vert^{2})^{2}}= \frac{J}{(1+\vert
\nabla u \vert^{2})^{2}}, $$
and the induced area form is
$$  dA= (1+\vert \nabla u \vert^{2})^{1/2} dx_{1} dx_{2}. $$
Thus
$$   K^{1/4} dA = J^{1/4} dx_{1} dx_{2}. $$
The left hand side is invariant under Euclidean transformations of $\bR^{3}$
while the right hand side is invariant under unimodular affine transformations
of $\bR^{2}$. Since these two groups generate all unimodular affine transformations
of $\bR^{3}$ we see that this $2$-form is an {\it affine} invariant of the surface.
The graphs of the solutions of  equation (2) when $\psi(t)=t^{1/4}$ are  \lq\lq affine
maximal'' surfaces in $\bR^{3}$. According to Chern and Terng \cite{kn:CT}
these surfaces can be described locally as follows. Let $F_{1}, F_{2}, F_{3}$ be harmonic
functions of variables $\lambda_{1}, \lambda_{2}$ and write $\uF=\uF(\lambda_{1},
\lambda_{2})$ for the
corresponding vector-valued function. Then the condition $\frac{\partial^{2}
\uF}{\partial \lambda_{1}^{2}}+\frac{\partial^{2} \uF}{\partial \lambda_{2}^{2}}=0$
implies the consistency of the first order system
\begin{equation} \frac{\partial \uZ}{\partial \lambda_{1}} = \uF \times \frac{\partial
\uF}{\partial \lambda_{2}} \ \ \ \ \ \frac{\partial \uZ}{\partial \lambda_{2}} =- \uF \times \frac{\partial
\uF}{\partial \lambda_{1}} \end{equation}
so there is a solution $\uZ$ which, under appropriate non-degeneracy conditions,
parametrises a surface in $\bR^{3}$. Chern and Terng show that these surfaces
are precisely the affine maximal surfaces. We want to relate this description
to ours. Notice first that our equations can be written in a similar form.
Given a pair of solutions $\xi_{1}, \xi_{2}$ of (3) we define a vector valued
function $\Xi$ with components $\xi_{1}, \xi_{2}, 1$. Then if $\uZ=(x_{1},
x_{2}, u)$ our description is the first order system
\begin{equation}   \frac{\partial \uZ}{\partial H}= p(r)\  \Xi\times \frac{\partial
\xi}{\partial r} \ \  \  \frac{\partial \uZ}{\partial r}= -p(r)\  \Xi\times \frac{\partial
\xi}{\partial H}. \end{equation}
In the case when $\psi(t)=t^{1/4}$ we get $p(r)=r^{2}$ so the equation (3)
is not  the ordinary Laplace equation, in the variables $(H,r)$. However
given a  function $F(H,r)$ we consider the first order system
\begin{equation}  p(r) \frac{\partial \xi}{\partial H}= r \frac{\partial F}{\partial H}
\  \ \ \ \ p(r) \frac{\partial \xi}{\partial r}= r \frac{\partial F}{\partial
r}- F. \end{equation}
If $p(r)=r^{2}$ one readily checks that this system is consistent, so there
is a solution $\xi(H,r)$. The equation (3) for $\xi$ is equivalent to 
the ordinary Laplace equation for $F$. Now take two harmonic functions $F_{1},
F_{2}$ and set $F_{3}(H,r)=r$. Then a few lines of calculation show that
the system (12) is identical to the system (13), when in the latter we use the
functions $\xi_{1}, \xi_{2} $ corresponding to $F_{1}, F_{2}$ by solving
(14).

\end{enumerate}


\end{document}